\documentclass[a4paper]{amsart}
\usepackage{amsmath,amssymb}

\numberwithin{equation}{section}

\newcommand{\nc}{\newcommand}
\nc{\C}{\mathbf{C}}		% complex number
\nc{\Z}{\mathbf{Z}}		% integer
\nc{\apo}{S}                    % antipode
\nc{\cop}{\Delta}               % coproduct
\nc{\Heis}{\mathcal{B}}		% Heisenberg algebra
\nc{\calV}{\mathcal{V}}
\nc{\half}{\frac{1}{2}}
\nc{\la}{\lambda}
\nc{\La}{\Lambda}
\nc{\ot}{\otimes}
\nc{\ra}{\rightarrow}
\nc{\uqp}{U'_q}
\nc{\uq}{U_q}
\nc{\slth}{\widehat{\mathfrak{sl}}_2}
\nc{\spfh}{\widehat{\mathfrak{sp}}_4}
\nc{\flip}{\sigma}		% Dynkin diagram automorphism
\nc{\ket}[1]{|#1\rangle}
\nc{\Philzz}{\Phi^{(\ell)}_{0}(z)}
\nc{\Philzw}{\Phi^{(\ell)}_{0}(w)}
\nc{\Psillz}{\Psi^{(\ell)}_{\ell}(z)}
\nc{\Psillw}{\Psi^{(\ell)}_{\ell}(w)}
\nc{\Phitzz}{\Phi^{(2)}_{0}(z)}
\nc{\Phitzw}{\Phi^{(2)}_{0}(w)}
\nc{\Psittz}{\Psi^{(2)}_{2}(z)}
\nc{\Psittw}{\Psi^{(2)}_{2}(w)}
\nc{\Dh}{D^{\frac{1}{2}}}
\nc{\Dmh}{D^{-\frac{1}{2}}}
\nc{\qexp}{\exp_q}
\nc{\qiexp}{\exp_{q^{-1}}}
\nc{\infp}[1]{(#1; q^4)_\infty}
\DeclareMathOperator{\End}{End}
\DeclareMathOperator{\Sym}{Sym}

\theoremstyle{plain}
\newtheorem{thm}{Theorem}[section]
\newtheorem{prop}[thm]{Proposition}
\newtheorem{lem}[thm]{Lemma}
\theoremstyle{definition}
\newtheorem{defi}[thm]{Definition}

\begin{document}

%%%%%%%%%%%%%%%%%%%%%%%%%%%%%%%%%%%%%%%%%%%%%%%%%%%%%%%%%%%%%%%%%%%%%%

\title[basic representations of $\uq(\spfh)$]
{A construction of\\
level $1$ irreducible modules for $\uq(\spfh)$\\
using level $2$ intertwiners for $\uq(\slth)$}

\author{Boris Feigin, Jin Hong, and Tetsuji Miwa}
\address{L. D. Landau Institute of Theoretical Physics\\
         Chernogolovka 142432\\
         Russian Federation}
\address{Korea Institute for Advanced Study\\
         207-43 Cheongryangri-dong, Dongdaemun-gu\\
         Seoul 130-012, Korea}
\address{Department of Mathematics\\
         Kyoto University\\
         Graduate School of Science\\
         Sakyo, Kyoto\\
         Japan}
%\thanks{}
\email{feigin@landau.ac.ru, jinhong@kias.re.kr, tetsuji@kusm.kyoto-u.ac.jp}

\begin{abstract}
We bosonize certain components of level $\ell$
$\uq(\slth)$-intertwiners of $(\ell + 1)$-dimensions.
For $\ell = 2$,
these intertwiners,
after certain modification by bosonic vertex operators,
are added to the algebra $\uq(\slth)$ at level $2$
to construct all irreducible highest weight representations of level $1$
for the quantum affine algebra $\uq(\spfh)$.
\end{abstract}

\maketitle

%%%%%%%%%%%%%%%%%%%%%%%%%%%%%%%%%%%%%%%%%%%%%%%%%%%%%%%%%%%%%%%%%%%%%%
\section{Introduction}

The aim of this paper is to construct the level $1$ irreducible
highest weight representations of $\uq(\spfh)$.
We start with the level $2$ irreducible highest weight
representations of $\uq(\slth)$.
We construct
a representation space of $\uq(\spfh)$ as an infinite sum of
representation spaces of $\uq(\slth)$ tensored by bosonic Fock
spaces.

We consider two sets of bosonic oscillators.
The first one is a part of the Drinfel$'\negthinspace$d
generators of $\uq(\slth)$ and the second one
is added by hand to extend the representation space as mentioned above.

On the extended space thus constructed, we define actions
of the Drinfel$'\negthinspace$d
generators of $\uq(\spfh)$.
{}For this purpose we use
the action of $\uq(\slth)$ as a part.
We also exploit the 
level $2$ intertwiners of dimensions $3$ in the sense
of~\cite{MR94c:17024,MR96e:82037}.
There are two types of such intertwiners, type I and II.
One of three components in each case can be written as
a vertex operator in terms of bosonic oscillators in $\uq(\slth)$.
We modify these operators by multiplying simple vertex operators
constructed from the second set of bosonic oscillators, and add the modified
vertex operators to the action.

The construction for the affine Lie algebra case, i.e., $q=1$,
was mentioned in~\cite{MR2000g:17042}.
Our construction is a $q$ deformation of their
construction. However, there are two technical differences.
First, in the $q$-deformed situation, the
type I and II intertwiners are not the same.
In our construction
we need both in a proper combination.
Second, in the $q$-deformed situation,
we construct only Drinfel$'\negthinspace$d
generators instead of constructing all the
$\mathfrak{sp}_4$ currents as in the affine Lie algebra case.
The cost is to prove
Drinfel$'\negthinspace$d $q$-Serre relations.

In~\cite{MR2001e:17024}, bosonizations of level $1$ representations
of $\uq(\widehat{\mathfrak{sp}}_{2n})$
were constructed.
Our construction is different
from theirs since we have constructed irreducible representations.

In Section~\ref{sec:2}, we prepare basic definitions for
$\uq(\slth)$ and $\uq(\spfh)$.
In Section~\ref{sec:3},
we construct special components of the type I and II intertwiners
for $\uq(\slth)$ at level $l$, in general.
In Section~\ref{sec:4},
we combine these to construct level $1$ irreducible
representations of $\uq(\spfh)$.

%%%%%%%%%%%%%%%%%%%%%%%%%%%%%%%%%%%%%%%%%%%%%%%%%%%%%%%%%%%%%%%%%%%%%%
\section{Quantum affine algebras and intertwiners}\label{sec:2}

Basic notations used in this paper is given in this section.
We shall deal with two quantum affine algebras,
$\uq(\slth)$ and $\uq(\spfh)$,
in this paper.
The quantum affine algebras $\uq$ of type $A_1^{(1)}$ and $C_2^{(1)}$
are generated by the elements
\begin{equation*}
e_i,\ f_i,\ t_i^{\pm1},\ q^{\pm d},
\end{equation*}
with $i=0,1$ for type $A_1^{(1)}$ and $i=0,1,2$ for type $C_2^{(1)}$.
We shall not write down their defining relations
in terms of these generators.
The subalgebras of $\uq$ generated by elements above, excluding the
elements $q^{\pm d}$,
will be denoted by $\uqp$.
Hopf algebra structure for both of these algebras is given as follows.
\begin{align*}
\cop(e_i) &= e_i\ot t_i + 1\ot e_i,\\
\cop(f_i) &= f_i\ot 1 + t_i^{-1}\ot f_i,\\
\cop(t_i) &= t_i\ot t_i,\\
\cop(q^d) &= q^d\ot q^d.
\end{align*}

The canonical central element for $\slth$ is given by
$c = h_0 + h_1$ and that of $\spfh$ is given by
$c = h_0 + h_1 + h_2$.

%%%%%%%%%%%%%%%%%%%%%%%%%%%%%%%%%%%%%%%%
\subsection{Intertwiners for $\uqp(\slth)$}

We now state the Drinfel$'\negthinspace$d
realization for the quantum affine algebra $\uq(\slth)$.
The generators are given by
\begin{equation*}
x^{\pm}(k),\ a(l),\ K^{\pm 1},\ \gamma^{\pm\half},\ q^{\pm d}.
\end{equation*}
Here, the indices run over $k\in\Z$ and $l\in\Z^{\times}$.
In terms of the generating function
\begin{equation*}
X^\pm(z) = \sum_{k\in\Z} x^\pm(k)\, z^{-k-1},
\end{equation*}
the defining relations are given by
{\allowdisplaybreaks
\begin{align}\text{}
&[\gamma^{\pm\half}, u] = 0 \quad\text{for all $u\in\uq(\slth)$},\\
&[a(k), a(l)] = \delta_{k+l,0} \frac{[2k]}{k}
  \frac{\gamma^k-\gamma^{-k}}{q-q^{-1}},\\
&K a(k) K^{-1} = a(k),\quad
  K X^\pm(z) K^{-1} = q^{\pm2} X^\pm(z),\\
&[a(k), X^{\pm}(z)] = \pm \frac{[2k]}{k} \gamma^{\mp|k|/2} z^k X^\pm(z),%
\label{eq:24}\\
&(z-q^{\pm2}w) X^\pm(z)X^\pm(w) + (w-q^{\pm2}z)X^\pm(w)X^\pm(z) = 0,%
\label{eq:25}\\
&[X^+(z),\;X^-(w)] \notag\\
&\phantom{X^+(z),\;X}
 = K \exp\bigg\{ (q-q^{-1})\sum_{k=1}^\infty a(k)\gamma^{k/2} z^{-k}\bigg\}
  \frac{\delta(z/\gamma w)}{(q-q^{-1})zw}\label{eq:26}\\
&\phantom{X^+(z),\;X}
 - K^{-1} \exp\bigg\{ -(q-q^{-1})\sum_{k=1}^\infty a(-k)\gamma^{k/2}z^k\bigg\}
  \frac{\delta(\gamma z/w)}{(q-q^{-1})zw},\notag\\
&q^d K q^{-d} = K,\ 
q^d x^{\pm}(k) q^{-d} = q^k x^\pm(k),\ 
 q^d a(k) q^{-d} = q^k a(k).
\end{align}
} % end of \allowdisplaybreaks
Here, the notation $\delta(z) = \sum_{k\in\Z} z^k$ is a formal infinite sum.
Identification between the two presentations of $\uq(\slth)$ is given
in~\cite{MR95i:17011,MR88j:17020,MR99j:17021}.
\begin{alignat*}{3}
K&\leftrightarrow t_1,
  &\quad \gamma &\leftrightarrow t_0 t_1,
  &\quad q^d &\leftrightarrow q^d,\\
x^+(0) &\leftrightarrow e_1,
  &&
  & x^-(0) &\leftrightarrow f_1,\\
x^+(-1) &\leftrightarrow t_0 f_0,
  &&
  &\quad x^-(1) &\leftrightarrow e_0 t_0^{-1}.
\end{alignat*}

The algebra $\uqp(\slth)$,
the subalgebra of $\uq(\slth)$ without $q^d$,
admits a finite dimensional representation.
\begin{equation*}
V_z^{(\ell)} = \bigoplus_{j=0}^\ell \C(q) v_j^{(\ell)}.
\end{equation*}
The action of $\uqp(\slth)$ on $V_z^{(\ell)}$ is given by
\begin{align*}
f_1 v_j^{(\ell)} &= [j+1] v_{j+1}^{(\ell)},\\
e_1 v_j^{(\ell)} &= [\ell-j+1] v_{j-1}^{(\ell)},\\
t_1 v_j^{(\ell)} &= q^{\ell-2j} v_j^{(\ell)},\\
f_0 v_j^{(\ell)} &= z^{-1}[\ell-j+1] v_{j-1}^{(\ell)},\\
e_0 v_j^{(\ell)} &= z\,[j+1] v_{j+1}^{(\ell)},\\
t_0 v_j^{(\ell)} &= q^{2j-\ell} v_j^{(\ell)}.
\end{align*}

Given a dominant integral weight $\la$ of level $\ell$, i.e., satisfying
$\la(h_0 + h_1) = \ell$, we let
\begin{align*}
\Phi^{(\ell)}(z)
   :\;& V^{(\ell)}_z\ot V(\la) \ra V(\flip\la), \quad\text{(type-I)}\\
\Psi^{(\ell)}(z)
   :\;& V(\la)\ot V^{(\ell)}_z \ra V(\flip\la), \quad\text{(type-II)}
\end{align*}
denote the intertwiners.
Here, the map $\flip$ permutes the fundamental weights $\La_0$ and $\La_1$
of $\uqp(\slth)$.
Note that these intertwiners are unique up to scalar multiple.
We define the components of the intertwiners by
\begin{align*}
\Phi_j^{(\ell)}(z) \ket{v} &= \Phi^{(\ell)}(z) v^{(\ell)}_j \ot\ket{v},\\
\Psi_j^{(\ell)}(z) \ket{v} &= \Psi^{(\ell)}(z) \ket{v}\ot v^{(\ell)}_j.
\end{align*}

%%%%%%%%%%%%%%%%%%%%%%%%%%%%%%%%%%%%%%%%
\subsection{Drinfel$'\negthinspace$d
realization for $\uq(\spfh)$}\label{ssec:22}

Let $A = (a_{i,j})_{i\in \hat{I}}$ be the Cartan matrix of type $C_2^{(1)}$.
Here, the index set $\hat{I} = \{0,1,2\}$.
We set $q_1 = q$, $q_2 = q^2$ and let $[n]_i$ denote the $q$-integer
which uses $q_i$ in place of $q$.

The Drinfel$'\negthinspace$d
generators for $\uq(\spfh)$ will be denoted by
\begin{equation*}
x_i^\pm(k),\ a_i(l),\ K_i^{\pm 1},\ \gamma^{\pm\half},\ q^{\pm d}.
\end{equation*}
Here, the indices run over $i\in I = \{1,2\}$, $k\in\Z$,
and $l\in\Z^{\times}$.
In terms of the generating function
\begin{equation*}
X_i^\pm(z) = \sum_{k\in\Z} x_i^\pm(k)\, z^{-k-1},
\end{equation*}
the defining relations are given by
{\allowdisplaybreaks
\begin{align}\text{}
&[\gamma^{\pm\half}, u]=0 \quad\text{for all $u\in\uq(\spfh)$},\\
&K_i K_j = K_j K_i, \quad K_i K_i^{-1} = K_i^{-1} K_i = 1,\\
&[a_i(k), a_j(l)] = \delta_{k+l,0} \frac{[a_{i,j}k]_i}{k}
  \frac{\gamma^k-\gamma^{-k}}{q_j-q_j^{-1}},\\
&K_j a_i(k) K_j^{-1} = a_i(k),\quad
  K_j X_i^\pm(z) K_j^{-1} = q_i^{\pm a_{i,j}} X_i^\pm(z),\\
&q^d K_i q^{-d} = K_i,\ 
  q^d x_i^{\pm}(k) q^{-d} = q^k x_i^\pm(k),\ 
  q^d a_i(k) q^{-d} = q^k a_i(k),\label{eq:n212}\\
&[a_i(k), X_j^{\pm}(z)] = \pm \frac{[a_{i,j}k]_i}{k}
  \gamma^{\mp|k|/2} z^k X_j^\pm(z),\label{eq:211}\\
&(z-q_i^{\pm a_{i,j}}w) X_i^\pm(z)X_j^\pm(w)
  + (w-q_i^{\pm a_{i,j}}z)X_j^\pm(w)X_i^\pm(z) = 0,\label{eq:212}\\
&[X_i^+(z), X_j^-(w)]
 = \frac{\delta_{i,j}}{q_i-q_i^{-1}} \notag\\
&\phantom{[X_i^+(z), X}
  \times\bigg(
  K_i \exp\bigg\{ (q_i-q_i^{-1})
  \sum_{k=1}^\infty a_i(k)\gamma^{k/2} z^{-k}\bigg\}
  \frac{\delta(z/\gamma w)}{zw}\label{eq:213}\\
&\phantom{[X_i^+(z), X}
- K_i^{-1} \exp\bigg\{ -(q_i-q_i^{-1})
  \sum_{k=1}^\infty a_i(-k)\gamma^{k/2}z^k\bigg\}
  \frac{\delta(\gamma z/w)}{zw}
  \bigg),\notag\\
&\Sym_{z_1,z_2,z_3}
\begin{pmatrix}\label{eq:214}
\hfill
        X_2^\pm(w) X_1^\pm(z_1) X_1^\pm(z_2) X_1^\pm(z_3)\\[0.1em]
\hfill
- [3]_1 X_1^\pm(z_1) X_2^\pm(w) X_1^\pm(z_2) X_1^\pm(z_3)\\[0.1em]
\hfill
+ [3]_1 X_1^\pm(z_1) X_1^\pm(z_2) X_2^\pm(w) X_1^\pm(z_3)\\[0.1em]
\hfill
      - X_1^\pm(z_1) X_1^\pm(z_2) X_1^\pm(z_3) X_2^\pm(w)
\end{pmatrix}
= 0,\\
&\Sym_{z_1,z_2}
\begin{pmatrix}\label{eq:215}
\hfill
        X_1^\pm(w) X_2^\pm(z_1) X_2^\pm(z_2)\\[0.1em]
\hfill
- [2]_2 X_2^\pm(z_1) X_1^\pm(w) X_2^\pm(z_2)\\[0.1em]
\hfill
      + X_2^\pm(z_1) X_2^\pm(z_2) X_1^\pm(w)
\end{pmatrix}
= 0.
\end{align}
} % end of \allowdisplaybreaks

We remark that identification of this realization
with the Drinfel$'\negthinspace$d
realization given in~\cite{MR2001e:17024} may be done by mapping
\begin{align*}
q^{\half} &\leftrightarrow q,\\
a_{1,k} &\leftrightarrow \frac{1}{[2]} a_1(k),\\
a_{2,k} &\leftrightarrow a_2(k),
\end{align*}
and mapping all other elements trivially.
Here, the left hand side denotes elements from~\cite{MR2001e:17024}
and the right hand side denotes elements from this paper.

%%%%%%%%%%%%%%%%%%%%%%%%%%%%%%%%%%%%%%%%%%%%%%%%%%%%%%%%%%%%%%%%%%%%%%
\section{Bosonization of intertwiners}\label{sec:3}

We start this section by
first recalling some facts from~\cite{MR95e:17021}.
Define endomorphisms $S_i$ ($i=0,1$)
acting on integrable $\uq(\slth)$-modules by
\begin{equation}
S_i = \qiexp(q^{-1}e_i t_i^{-1})\qiexp(-f_i)\qiexp(q e_i t_i)
      q^{h_i(h_i +1)/2}.
\end{equation}
Here, we have used the formal notation
\begin{equation*}
\qexp(x) = \sum_{n=0}^{\infty} \frac{q^{\half n(n-1)}}{[n]_{q}!} x^n.
\end{equation*}
We remark that
$\qexp(x) \qiexp(-x) = 1$
so that $S_i$ is invertible.

\begin{prop}[\cite{MR95e:17021}]\label{pr:31}
The endomorphism $S_i$ satisfies the following relations.
\begin{align}
S_i e_i S_i^{-1} &= - f_i t_i,\label{eq:n32}\\
S_i f_i S_i^{-1} &= - t_i^{-1} e_i,\label{eq:n33}\\
S_i t_i S_i^{-1} &= t_i^{-1},\label{eq:n34}\\
S_i e_j S_i^{-1} &= \frac{1}{[2]}
  ( q^{-2} e_j e_i^2 - q^{-1} [2] e_i e_j e_i + e_i^2 e_j ),\label{eq:n35}\\
S_i f_j S_i^{-1} &= \frac{1}{[2]}
  ( q^2 f_i^2 f_j - q [2] f_i f_j f_i + f_j f_i^2 ),\\
S_i t_j S_i^{-1} &= t_j t_i^{2}.
\end{align}
Here, we only take $i\neq j$.
\end{prop}

%%%%%%%%%%%%%%%%%%%%%%%%%%%%%%%%%%%%%%%%
\subsection{The operator $\Dh$}

We now set
\begin{equation}
\Dh = S_0 t_0^{-1} \flip.
\end{equation}
Recall that $\flip$ permutes the fundamental weights $\La_0$ and $\La_1$.
Here, the map
\begin{equation*}
\flip : V(\la) \ra V(\sigma \la)
\end{equation*}
sends $v_{\la}$ to $v_{\sigma\la}$ and
satisfies $\sigma f_i = f_{1-i} \sigma$.
Hence the operator $\Dh$ is a map from $V(\la)$ to $V(\flip\la)$.

\begin{prop}
On the Drinfel\/$'\negthinspace$d
generators, the operator $\Dh$ has the following properties.
\begin{align}
\Dh a(n) \Dmh &= a(n),\label{eq:n309}\\
\Dh K \Dmh &= \gamma^{-1} K,\label{eq:n310}\\
\Dh X^{\pm}(z) \Dmh &= -z^{\mp 1} X^{\pm}(z).\label{eq:n311}
\end{align}
\end{prop}
\begin{proof}
Equation~\eqref{eq:n310} is immediate from~\eqref{eq:n34}.

Let us show the last one.
The two special cases
\begin{align}
\Dh x^+(0) \Dmh &= - x^+(-1),\label{eq:n312}\\
\Dh x^-(0) \Dmh &= - x^-(1)\label{eq:n313t}
\end{align}
may be obtained from~\eqref{eq:n32} and~\eqref{eq:n33}.
 From the defining relations for $\uq(\slth)$, we may write
\begin{equation*}
x^-(2)
 = - \frac{\gamma^{-\half}}{[2]}[ a(1) , x^-(1)]
 = - \frac{\gamma^{-1}}{[2]} [ (e_1 e_0 - q^2 e_0 e_1), e_0 t_0^{-1}].
\end{equation*}
This looks quite similar to the right hand side of~\eqref{eq:n35}.
Indeed, if we compare this with the outcome of applying
Proposition~\ref{pr:31} to $\Dh X^-(1) \Dmh$, we obtain
\begin{equation}
\Dh x^-(1) \Dmh = - x^-(2).\label{eq:n313}
\end{equation}
Equations~\eqref{eq:n312} and~\eqref{eq:n313} may then be applied to
\begin{align*}
[ x^+(0), x^-(1) ] = \gamma^{-\half} K a(1)
\end{align*}
to show
\begin{equation}
\Dh a(1) \Dmh = a(1).
\end{equation}
With this, starting from~\eqref{eq:n312} and~\eqref{eq:n313t},
we may recursively show~\eqref{eq:n311} true for all
negative powers of $z$.
Positive powers of~\eqref{eq:n311} may similarly be shown
after obtaining $\Dh a(-1) \Dmh = a(-1)$.

It remains to deal with the first equation.
To this end, recall that one of the defining relations~\eqref{eq:26}
for $\uq(\slth)$ may be written as
\begin{equation}\label{eq:n316}
[x^+(k), x^-(l)] = 
 \frac{\gamma^{\frac{k-l}{2}}\psi_{k+l} - \gamma^{\frac{l-k}{2}}\varphi_{k+l}}
      {q - q^{-1}},
\end{equation}
where
\begin{align*}
\sum_{k=0}^\infty \psi_k z^{-k}
   &= K \exp\bigg\{(q-q^{-1})\sum_{k=1}^\infty a(k)\, z^{-k} \bigg\},\\
\sum_{k=0}^\infty \varphi_{-k} z^k
   &= K^{-1} \exp\bigg\{ -(q-q^{-1})\sum_{k=1}^\infty a(-k)\, z^k \bigg\}.
\end{align*}
Applying~\eqref{eq:n311} to~\eqref{eq:n316}, we may show
$\Dh \psi_k \Dmh = \gamma^{-1}\psi_k$ for all nonnegative $k$.
Since the $\gamma^{-1}$ term comes from the $K$ inside $\psi_k$,
we may start with
$\Dh a(1) \Dmh = a(1)$, which we already obtained, to recursively show
\begin{equation*}
\Dh a(k) \Dmh = a(k)
\end{equation*}
for all positive $k$.
Results for negative $k$ may be obtained similarly.
\end{proof}

%%%%%%%%%%%%%%%%%%%%%%%%%%%%%%%%%%%%%%%%
\subsection{Bosonization of $\Phi_0^{(\ell)}(z)$.}

We shall follow~\cite{MR90e:17028,MR96e:82037} in realizing the
highest component $\Phi_0^{(\ell)}$ of the type-I $\uqp(\slth)$-intertwiner.

Since all components of the intertwiner act on a level $\ell$
module, we will freely use $\gamma = q^\ell$.
Starting from the fact that $\Phi^{(\ell)}(z)$ commutes with the
$\uqp(\slth)$-action, we obtain the following identities.
\begin{align}
[ \ell-j+1 ]\; \Phi^{(\ell)}_{j-1}(z)\; t_1
  &= e_1 \Phi^{(\ell)}_{j}(z) - \Phi^{(\ell)}_{j}(z)e_1,\label{eq:31}\\
z\,[ j+1 ]\; \Phi^{(\ell)}_{j+1}(z)\; t_0
  &= e_0 \Phi^{(\ell)}_{j}(z) - \Phi^{(\ell)}_{j}(z)e_0,\label{eq:32}\\
[ j+1 ]\; \Phi^{(\ell)}_{j+1}(z)
  &= f_1 \Phi^{(\ell)}_{j}(z) - q^{(2j-\ell)}\Phi^{(\ell)}_{j}(z)f_1,%
\label{eq:33}\\
z^{-1} [ \ell-j+1 ]\; \Phi^{(\ell)}_{j-1}(z)
  &= f_0 \Phi^{(\ell)}_{j}(z)
     - q^{(\ell-2j)}\Phi^{(\ell)}_{j}(z)f_0,\label{eq:34}\\
t_1 \Phi^{(\ell)}_{j}(z) t_1^{-1} 
  &= q^{(\ell-2j)} \Phi^{(\ell)}_{j}(z),\label{eq:35}\\
t_0 \Phi^{(\ell)}_{j}(z) t_0^{-1} 
  &= q^{(2j-\ell)} \Phi^{(\ell)}_{j}(z).\label{eq:36}
\end{align}

 From~\eqref{eq:31} and~\eqref{eq:34} we get
\begin{align}
[ x^+(0), \Philzz ] &= 0,\label{eq:37}\\
[ x^+(-1), \Philzz ] &= 0.\label{eq:38}
\end{align}
We may combine the equations~\eqref{eq:32} and~\eqref{eq:33}
by removing the left hand side.
After changing everything into the Drinfel$'\negthinspace$d notation we get
\begin{equation}\label{eq:310}
z\Big( \Philzz\,x^-(0) - \gamma\; x^-(0)\,\Philzz \Big)
= \Big( \gamma\;\Philzz\, x^-(1) - x^-(1)\,\Philzz \Big).
\end{equation}
With equations~\eqref{eq:37}, \eqref{eq:38}, and \eqref{eq:310}, we may
proceed as in~\cite[pp. 73--74]{MR96e:82037} to obtain the following
proposition
\begin{prop}\label{pr:33}
The component $\Philzz$ of the type-I intertwiner satisfies
the following set of equations.
\begin{align}
&K \Philzz K^{-1} = \gamma \Philzz,\\
&[ X^+(w), \Philzz ] = 0,\label{eq:327}\\
&[ a(\pm k), \Philzz ]
   = \frac{1}{k} \gamma^{\frac{k}{2}} [ k\ell ]\, z^{\pm k} \Philzz.
\end{align}
\end{prop}

This proposition allows us to guess a realization for $\Philzz$.

\begin{thm}\label{thm:34}
The component $\Philzz$ of the type-I intertwiner may be
realized, up to scalar multiple, as follows.
\begin{align*}
\Philzz =&
   \exp\bigg( \sum_{k=1}^{\infty}
     \frac{\gamma^{\frac{k}{2}}}{[2k]} a(-k)\, z^k \bigg)
   \exp\bigg( - \sum_{k=1}^{\infty}
     \frac{\gamma^{\frac{k}{2}}}{[2k]} a(k)\, z^{-k} \bigg)\\
   &\times\Dh z^{\half h_1}.
\end{align*}
\end{thm}
\begin{proof}
That this definition of $\Philzz$
satisfies the three equations of Proposition~\ref{pr:33} may be
checked through routine calculations.

To prove that this is the correct realization for $\Philzz$,
first recall that the group of equations~\eqref{eq:31}--\eqref{eq:36}
is equivalent to the definition of the intertwiner $\Phi^{(\ell)}(z)$.
So let us define all other components of the intertwiner
by~\eqref{eq:33} and show that the defined components
satisfy the equations~\eqref{eq:31}--\eqref{eq:36}.

Checking equations~\eqref{eq:35} and~\eqref{eq:36} is trivial.
Let us consider~\eqref{eq:31}.
The $j=0$ case is immediate from~\eqref{eq:327}.
So let us assume~\eqref{eq:31} true for $j=s$.
Then, using $j=s$ and $j=s-1$ cases of~\eqref{eq:33},
we have
\begin{align*}
[s&+1] ( e_1 \Phi^{(\ell)}_{s+1}(z) - \Phi^{(\ell)}_{s+1}(z) e_1 )\\
&= e_1 f_1 \Phi^{(\ell)}_s(z) - q^{(2s-\ell)} e_1\Phi^{(\ell)}_s(z) f_1
   - f_1 \Phi^{(\ell)}_s(z) e_1 + q^{(2s-\ell)} \Phi^{(\ell)}_s(z) f_1 e_1\\
&= ( [e_1, f_1]\,\Phi^{(\ell)}_{s}(z)
     - q^{(2s-\ell)}\Phi^{(\ell)}_{s}(z) [e_1, f_1] )\\
&\phantom{==}+ [\ell-s+1]( f_1\Phi^{(\ell)}_{s-1}(z)t_1
                 - q^{(2s-\ell)}\Phi^{(\ell)}_{s-1}(z)t_1 f_1 )\\
&= [\ell-2s]\,\Phi^{(\ell)}_{s}(z) t_1
   + [\ell-s+1][s]\,\Phi^{(\ell)}_{s}(z) t_1\\
&= [s+1][\ell-s]\,\Phi^{(\ell)}_s t_1.
\end{align*}
This shows the induction step and~\eqref{eq:31} is true for all $j$.

In a very similar way, we may show~\eqref{eq:34} true for all $j$,
if we can show~\eqref{eq:32} true for all $j$.

It only remains to show~\eqref{eq:32}.
We prove it by induction on $j$.
The $j=0$ case of~\eqref{eq:32} is equivalent to the validity
of~\eqref{eq:310}.
And the induction step for~\eqref{eq:32} may be shown using~\eqref{eq:33}
as in the previous two induction proofs.
So this last step reduces to showing~\eqref{eq:310}.
Using~\eqref{eq:24} on
the explicit formula for $\Philzz$, given in the statment of this theorem,
we may calculate
\begin{align*}
&(1 - \gamma w/z)
 \exp\bigg( - \sum_{k=1}^{\infty}
     \frac{\gamma^{\frac{k}{2}}}{[2k]} a(k)\, z^{-k} \bigg)
 X^-(w)\\
&\phantom{mmmmmmmmmm}=
 X^-(w)
 \exp\bigg( - \sum_{k=1}^{\infty}
     \frac{\gamma^{\frac{k}{2}}}{[2k]} a(k)\, z^{-k} \bigg),\\
&(1 - \gamma z/w)
 X^-(w)
 \exp\bigg( \sum_{k=1}^{\infty}
     \frac{\gamma^{\frac{k}{2}}}{[2k]} a(-k)\, z^k \bigg)\\
&\phantom{mmmmmmmmmm}=
\exp\bigg( \sum_{k=1}^{\infty}
     \frac{\gamma^{\frac{k}{2}}}{[2k]} a(-k)\, z^k \bigg)
X^-(w).
\end{align*}
With this, we can show
\begin{equation*}
\big({z-\gamma w}\big)\Philzz X^-(w) 
= - \big({w-\gamma z}\big) X^-(w)\Philzz.
\end{equation*}
It leads to
\begin{align*}
z\Big( \Philzz X^-(w) -& \gamma X^-(w)\Philzz \Big)\\
 &= w\Big( \gamma\Philzz X^-(w) - X^-(w)\Philzz \Big),
\end{align*}
which is a generalization of~\eqref{eq:310}.
\end{proof}

%%%%%%%%%%%%%%%%%%%%%%%%%%%%%%%%%%%%%%%%
\subsection{Bosonization of $\Psi^{(\ell)}_{\ell}(z)$.}

Realization for the lowest component $\Psillz$ of the type-II intertwiner
will be given in this subsection.
Since all the steps are as in the previous subsection, we shall be
very brief.

{\allowdisplaybreaks
\begin{align}
[\ell-j+1]\Psi^{(\ell)}_{j-1}(z)
   &= e_1 \Psi^{(\ell)}_j(z) - q^{(\ell-2j)}\Psi^{(\ell)}_j e_1,\\
z\,[j+1]\Psi^{(\ell)}_{j+1}(z)
   &= e_0 \Psi^{(\ell)}_j(z) - q^{(2j-\ell)}\Psi^{(\ell)}_j(z) e_0,\\
[j+1]\Psi^{(\ell)}_{j+1}(z) t_1^{-1}
   &= f_1 \Psi^{(\ell)}_j(z) - \Psi^{(\ell)}_j(z) f_1,\\
z^{-1}[\ell-j+1]\Psi^{(\ell)}_{j-1}(z) t_0^{-1}
   &= f_0 \Psi^{(\ell)}_j(z) - \Psi^{(\ell)}_j(z) f_0,\\
t_1 \Psi^{(\ell)}_j t_1^{-1}
   &= q^{(\ell-2j)}\Psi^{(\ell)}_j,\\
t_0 \Psi^{(\ell)}_j t_0^{-1}
   &= q^{(2j-\ell)}\Psi^{(\ell)}_j.
\end{align}
} % end of \allowdisplaybreaks.

\begin{prop}
The component $\Psillz$ of the type-II intertwiner satisfies
the following set of equations.
\begin{align}
&K \Psillz K^{-1} = \gamma^{-1} \Psillz,\\
&[ X^-(w), \Psillz ] = 0,\label{eq:336}\\
&[ a(\pm k), \Psillz ]
  = -\frac{1}{k} \gamma^{-\frac{k}{2}} [ k\ell ]\, (q^2 z)^{\pm k} \Psillz.
\end{align}
\end{prop}

\begin{thm}
The component $\Psillz$ of the type-II intertwiner may be
realized, up to scalar multiple, as follows.
\begin{align*}
\Psillz =&
   \exp\bigg( - \sum_{k=1}^{\infty}
     \frac{\gamma^{-\frac{k}{2}}}{[2k]} a(-k) (q^2 z)^k\bigg)
   \exp\bigg( \sum_{k=1}^{\infty}
     \frac{\gamma^{-\frac{k}{2}}}{[2k]} a(k) (q^2 z)^{-k} \bigg)\\
   &\times\Dmh (q^2 z)^{-\half h_1}.
\end{align*}
\end{thm}

%%%%%%%%%%%%%%%%%%%%%%%%%%%%%%%%%%%%%%%%
\subsection{Relations for vertex operators}

With the explicit bosonizations obtained in the previous sections,
we may write down relations for these vertex operators.

\begin{lem}\label{lem:n38}
We have the following commutation relation between vertex operators
and the grading operator.
\begin{align*}
q^d \Phitzz q^{-d}
&= q^{\half\la(h_0) - 1} \Phi^{(2)}_0(q^{-1}z),\\
q^d \Psittz q^{-d}
&= q^{-\half\la(h_1)} \Psi^{(2)}_2(q^{-1}z),
\end{align*}
as operators from $V(\la)$ to $V(\sigma\la)$.
\end{lem}
\begin{proof}
We shall deal with just the first one.
Second one may be done similarly.
Using the explicit realization of $\Phitzz$, we may
reduce this proof to showing
\begin{equation*}
q^d \Dh z^{\half h_1} q^{-d}
= q^{\half\la(h_0) - 1} \Dh (q^{-1}z)^{\half h_1}.
\end{equation*}
To show this,
we first follow the next set of equalities.
\begin{align*}
q^d \Dh z^{\half h_1} q^{-d}
&= q^d S_0 t_0^{-1} \sigma z^{\half h_1} q^{-d}\\
&= q^d \sigma q^{-d} S_1 t_1^{-1} z^{\half h_1}\\
&= q^d \sigma q^{-d} \sigma \Dh z^{\half h_1}.
\end{align*}
Now, we can show that,
as operators on an irreducible highest weight module of highest weight $\mu$,
\begin{equation*}
q^d \sigma q^{-d} \sigma =  q^{\half\mu(h_1)} q^{-\half h_1}.
\end{equation*}
It can be done by
checking the action of both sides on a weight vector of weight
$\mu - (x \alpha_0 + y\alpha_1)$.
Recalling that $\Dh$ is a map from $V(\la)$ to $V(\sigma\la)$,
we can continue as follows.
\begin{align*}
q^d \Dh z^{\half h_1} q^{-d}
&= q^{\half\sigma\la(h_1)} q^{-\half h_1} \Dh z^{\half h_1}\\
&= q^{\half\sigma\la(h_1)} q^{-1} \Dh q^{-\half h_1} z^{\half h_1}\\
&= q^{\half\la(h_0) - 1} \Dh (q^{-1}z)^{\half h_1}.
\end{align*}
The proof is complete.
\end{proof}

For later calculations, we need the following \emph{normal ordering}
symbols.
\begin{align*}
&: a(k) a(l) : \ = \
\begin{cases}
a(k) a(l) &\text{if $k<0$},\\
a(l) a(k) &\text{if $k>0$},
\end{cases}\\
&:\Dh h_1: \ = \ :h_1 \Dh: \ = \ \Dh h_1.
\end{align*}
For example, we have
\begin{equation*}
z^{\half h_1} \Dh = z : z^{\half h_1} \Dh :.
\end{equation*}
We also use the formal notation
\begin{equation*}
(x;p) = \prod_{k=0}^\infty (1 - x p^k).
\end{equation*}

\begin{lem}
We have the following identities concerning normal orderings.
\begin{align*}
\Philzz\Philzw
  &= z^{\half\ell}
     \frac{\infp{q^2 w/z}}{\infp{q^{2+2\ell} w/z}}
     :\Philzz\Philzw:\,,\\
\Philzz\Psillw
  &= \Big(\frac{1}{z}\Big)^{\half\ell}
     \frac{\infp{q^{4+\ell} w/z}}{\infp{q^{4-\ell} w/z}}
     :\Philzz\Psillw:\,,\\
\Psillz\Philzw
  &= \Big(\frac{1}{q^2 z}\Big)^{\half\ell}
     \frac{\infp{q^{\ell} w/z}}{\infp{q^{-\ell} w/z}}
     :\Psillz\Philzw:\,,\\
\Psillz\Psillw
  &= (q^2 z)^{\half\ell}
     \frac{\infp{q^{2-2\ell} w/z}}{\infp{q^{2} w/z}}
     :\Psillz\Psillw:.
\end{align*}
\end{lem}
\begin{proof}
This may be easily verified by applying
the following standard formula.
\begin{equation*}
\exp(A)\exp(B) = \exp([A,B]) \exp(B)\exp(A),
\quad\text{when $[A,B]$ is a scalar.}
\end{equation*}
In simplifying the outcome, the identity
\begin{equation*}
- \sum_{k=1}^\infty \frac{1}{k} \frac{[\ell k]}{[2k]} z^k
= \log \frac{\infp{q^{2-\ell} z}}{\infp{q^{2+\ell} z}}
\end{equation*}
will be useful.
\end{proof}

For later use, we write down the $\ell = 2$ case of this lemma.
\begin{lem}
We have the following identities concerning normal orderings
of level~$2$ vertex operators.
\begin{align*}
\Phitzz\Phitzw
  &= z\,(1 - q^2 w/z) :\Phitzz\Phitzw:\,,\\
\Phitzz\Psittw
  &= \frac{1}{z}\frac{1}{(1 - q^2 w/z)} :\Phitzz\Psittw:\,,\\
\Psittz\Phitzw
  &= \frac{1}{q^2 z}\frac{1}{(1 - w/q^2 z)} :\Psittz\Phitzw:\,,\\
\Psittz\Psittw
  &= q^2 z\,(1 - w/q^2z) :\Psittz\Psittw:.
\end{align*}
\end{lem}

%%%%%%%%%%%%%%%%%%%%%%%%%%%%%%%%%%%%%%%%%%%%%%%%%%%%%%%%%%%%%%%%%%%%%%
\section{Extended actions of vertex operators ($\ell=2$)}\label{sec:4}

In this section, we define an algebra $U'$.
This is done by adding two vertex operators
to the algebra $\uqp(\slth)$ after tensoring by some correction terms.
 From now on, we shall assume $\ell = 2$.
This is equivalent to saying $\gamma = q^2$ in $\uq(\slth)$.

%%%%%%%%%%%%%%%%%%%%%%%%%%%%%%%%%%%%%%%%
\subsection{Heisenberg algebra}

Let us define a Heisenberg algebra $\Heis'$.
It is to be generated by the element $b(k)$ with $k\in\Z$.
Defining relations are given below.
\begin{align}
&[b(k), b(l)] = \delta_{k+l,0}\frac{1}{k}(q^{2k} - 1 + q^{-2k})
\quad (k\neq 0),\\
&[b(0), b(k)] = 0.
\end{align}

We may consider the representation $\calV_{p}$,
over $\Heis'$, which contains
the vacuum vector $v(p)$ satisfying
\begin{align*}
b(k) v(p) &= 0 \quad\text{for $k>0$},\\
b(0) v(p) &= p v(p).
\end{align*}
The operator $T^r$ acts from $\calV_p$ to $\calV_{p+r}$.
It commutes with all $b(k)$ ($k\neq 0$) and sends
$v(p) \mapsto v(p+r)$.
We may trivially check
\begin{equation}
\text{}
[b(0), T^r] = r.
\end{equation}

We define two operators acting on $\bigoplus_{p\in\Z}\calV_p$.
\begin{align*}
\Omega_0(z) &=
\exp\bigg(\sum_{k=1}^\infty b(-k) q^{k} z^k\bigg)
\exp\bigg(-\sum_{k=1}^\infty b(k) q^{k} z^{-k}\bigg)
T^1 z^{b(0)},\\
\Omega_2(z) &=
\exp\bigg(-\sum_{k=1}^\infty b(-k) q^{-k} z^k\bigg)
\exp\bigg(\sum_{k=1}^\infty b(k) q^{-k} z^{-k}\bigg)
T^{-1} z^{-b(0)}.
\end{align*}

\emph{Normal ordering} in the algebra $\Heis'$ is defined by
\begin{align*}
&: b(k) b(l) : \ = \
\begin{cases}
b(k) b(l) &\text{if $k<0$},\\
b(l) b(k) &\text{if $k>0$},
\end{cases}\\
&:T^r b(0): \ = \ :b(0) T^r: \ = \ T^r b(0).
\end{align*}

\begin{lem}
Normal ordering of products of the operators $\Omega_0$ and $\Omega_2$
are as follows.
\begin{align}
\Omega_0(z) \Omega_0(w)
  &= z \frac{(1 - w/z)(1 - q^4 w/z)}{1 - q^2 w/z}
     :\Omega_0(z) \Omega_0(w):\,,\\
\Omega_0(z) \Omega_2(w)
  &= \frac{1}{z}\frac{1- w/z}{(1 - w/q^2 z)(1 - q^2 w/z)}
     :\Omega_0(z) \Omega_2(w):\,,\\
\Omega_2(z) \Omega_0(w)
  &= \frac{1}{z}\frac{1 - w/z}{(1 - w/q^2 z)(1 - q^2 w/z)}
     :\Omega_2(z) \Omega_0(w):\,,\\
\Omega_2(z) \Omega_2(w)
  &= z \frac{(1 - w/z)(1 - w/q^4 z)}{1 - w/q^2 z}
     :\Omega_2(z) \Omega_2(w):.
\end{align}
\end{lem}

%%%%%%%%%%%%%%%%%%%%%%%%%%%%%%%%%%%%%%%%
\subsection{The algebra $U'$}

We finally set
\begin{align}
Y^+(z) &= \Psi^{(2)}_2(q^{-2}z)\ot\Omega_2(z),\\
Y^-(z) &= \Phitzz\ot\Omega_0(z).
\end{align}
We write
\begin{equation*}
Y^\pm(z) = \sum_{k\in\Z} y^\pm (k)\, z^{-k-1}.
\end{equation*}
They are operators acting on the following spaces,
or on their direct sums.
\begin{align*}
\calV(0)
&= \Big( V(2\La_0) \otimes \bigoplus_{p\in 2\Z} \calV_p \Big)
   \oplus
   \Big( V(2\La_1) \otimes \bigoplus_{p\in 2\Z} \calV_{p+1} \Big)\\
\calV(1)
&= \Big( V(\La_0 + \La_1) \otimes \bigoplus_{p\in \Z} \calV_{p+\half} \Big)\\
\calV(2)
&= \Big( V(2\La_0) \otimes \bigoplus_{p\in 2\Z} \calV_{p+1} \Big)
   \oplus
   \Big( V(2\La_1) \otimes \bigoplus_{p\in 2\Z} \calV_{p} \Big)
\end{align*}
Note that each $\calV = \calV(j)$ ($j=0,1,2$) may be seen both as a
$\uqp(\slth)$-module and as a $\Heis'$-module.

\begin{defi}
The algebra $U' = U'(j)$ is defined to be
the subalgebra of $\End(\calV(j))$
generated by all elements of the quantum affine algebra $\uqp(\slth)$,
all elements of the Heisenberg algebra $\Heis'$, and
all coefficients of the modified vertex operators $Y^\pm(z)$.
\end{defi}

\begin{lem}\label{lem:this}
For each $p\in \Z$, the following equations are true up to nonzero
scalar multiple.
\begin{align*}
v_{2\La_1} \ot v(p+1)
&= y^-(-(p+1)) \left( v_{2\La_0} \ot v(p) \right),\\
v_{\La_0 + \La_1} \ot v(p + \textstyle{\half})
&= x^-(1)\, y^-(-(p+1))
   \left( v_{\La_0 + \La_1} \ot v(p - \textstyle{\half}) \right),\\
v_{2\La_0} \ot v(p+1)
&= (x^-(1))^2\, y^-(-(p+2)) \left( v_{2\La_1} \ot v(p) \right),\\
v_{2\La_1} \ot v(p-1)
&= (x^+(0))^2\, y^+(p-1) \left( v_{2\La_0} \ot v(p) \right),\\
v_{\La_0 + \La_1} \ot v(p - \textstyle{\half})
&= x^+(0)\, y^+(p)
   \left( v_{\La_0 + \La_1} \ot v(p + \textstyle{\half}) \right),\\
v_{2\La_0} \ot v(p-1)
&= y^+(p) \left( v_{2\La_1} \ot v(p) \right).
\end{align*}
\end{lem}
\begin{proof}
We shall prove just the second one.
Other cases are similar.

We may calculate
\begin{equation*}
y^-( -(p+1))\left( v_{\La_0 + \La_1} \ot v(p - \textstyle{\half}) \right)
= \Dh v_{\La_0 + \La_1} \ot v(p + \textstyle{\half}).
\end{equation*}
We claim that
\begin{equation*}
\Dh v_{\La_0 + \La_1} \in V(\La_0 + \La_1)_{\La_0 + \La_1 - \alpha_0}.
\end{equation*}
Note that the space on the right is of dimension $1$,
so that, if our claim is correct,
$\Dh v_{\La_0 + \La_1}$ is a nonzero scalar multiple of
$x^+(-1) v_{\La_0 + \La_1}$.
Also note that
\begin{equation*}
x^-(1) x^+(-1) v_{\La_0 + \La_1}
= v_{\La_0 + \La_1}.
\end{equation*}
So our claim is equivalent to the statement given in this Lemma.

The claim is proved by computing the weight of $\Dh v_{\La_0 + \La_1}$,
using equation~\eqref{eq:n310} and Lemma~\ref{lem:n38}.
\end{proof}

Recall that $V(\la)$ is irreducible under
the action of $\uqp(\slth)$ and
that each $\calV_p$ is also irreducible under the action of $\Heis'$.
So each of the tensored spaces $V(\la)\ot\calV_p$ are irreducible
under the action of $\uqp(\slth) \cup \Heis'$.
Now, Lemma~\ref{lem:this} shows that $Y^\pm(z)$ links these irreducible
tensored spaces.
Thus we have obtained the following proposition.

\begin{prop}
Each $\calV(j)$ \textup{(}$j=0,1,2$\textup{)}
is irreducible under the action of $U'$.
\end{prop}

For later use, we state some calculation results
that may be obtained from lemmas concerning normal ordering,
appearing in previous sections.

We shall write
\begin{align*}
&: a(k) x^\pm(l) : \ = \ :x^\pm(l) a(k): \ = \ 
\begin{cases}
a(k) x^\pm(l) &\text{if $k<0$},\\
x^\pm(l) a(k) &\text{if $k>0$},
\end{cases}\\
&:\Dh x^\pm(l): \ = \ :x^\pm(l) \Dh: \ = \ x^\pm(l)\Dh,\\
&:h_1 x^\pm(l): \ = \ :x^\pm(l) h_1: \ = \ x^\pm(l) h_1.
\end{align*}

{\allowdisplaybreaks
\begin{align}
&X^+(z)\,Y^+(w) = \frac{1}{1 - w/q^2 z} :X^+(z)\,Y^+(w):.\label{eq:410}\\
&Y^+(z)\,X^+(w) = - \frac{w}{z} \frac{1}{(1 - w/q^2 z)}:Y^+(z)\,X^+(w):.\\
&X^-(z)\,Y^-(w) = \frac{1}{1 - q^2 w/z} :X^-(z)\,Y^-(w):.\\
&Y^-(z)\,X^-(w) = - \frac{w}{z} \frac{1}{(1 - q^2w/z)} :Y^-(z)\,X^-(w):.\\
&Y^+(z)\,Y^+(w) = (z - q^{-4}w)(z - w) :Y^+(z)\,Y^+(w):.\label{eq:414}\\
&Y^-(z)\,Y^-(w) = (z - w)(z - q^{4}w) :Y^-(z)\,Y^-(w):.\label{eq:417}\\
&Y^+(z)\,Y^-(w) = \frac{1}{z^2}\frac{1}{(1 - q^2 w/z)(1 - w/q^2 z)}
                 :Y^+(z)\,Y^-(w):.\label{eq:415}\\
&Y^-(z)\,Y^+(w) = \frac{1}{z^2}\frac{1}{(1 - q^2 w/z)(1 - w/q^2 z)}
                 :Y^-(z)\,Y^+(w):.\label{eq:416}\\
&[ a(k), Y^\pm(z) ] = \mp \frac{[2k]}{k} q^{\mp |k|} z^k Y^\pm(z).%
\label{eq:418}\\
&[ b(k), Y^\pm(z) ] = \mp \frac{(q^{2k} - 1 + q^{-2k})}{|k|} q^{\mp |k|}
    z^k Y^\pm(z).\label{eq:419}\\
&X^\pm(w)\,Y^\pm(z_1)\,Y^\pm(z_2)\notag\\*
&\phantom{X^\pm(w)}
 =
 \frac{1}{(1 - q^{\mp 2}z_1/w)}\frac{1}{(1 - q^{\mp 2}z_2/w)}\label{eq:420}\\*
&\phantom{X^\pm(w) =\ }
\times (z_1 - q^{\pm 4}z_2)(z_1 - z_2)
   :X^\pm(w)\,Y^\pm(z_1)\,Y^\pm(z_2):.\notag\\
&Y^\pm(z_1)\,X^\pm(w)\,Y^\pm(z_2)\notag\\*
&\phantom{X^\pm(w)}
 = - \frac{w}{z_1}
 \frac{1}{(1 - q^{\mp 2}w/z_1)}\frac{1}{(1 - q^{\mp 2}z_2/w)}\label{eq:421}\\*
&\phantom{X^\pm(w) =\ }
   \times (z_1 - q^{\mp 4}z_2)(z_1 - z_2)
   :Y^\pm(z_1)\,X^\pm(w)\,Y^\pm(z_2):.\notag\\
&Y^\pm(z_1)\,Y^\pm(z_2)\,X^\pm(w)\notag\\*
&\phantom{X^\pm(w)}
 = \frac{w^2}{z_1 z_2}
 \frac{1}{(1 - q^{\mp 2}w/z_1)}\frac{1}{(1 - q^{\mp 2}w/z_2)}\label{eq:422}\\*
&\phantom{X^\pm(w) =\ }
   \times (z_1 - q^{\mp 4}z_2)(z_1 - z_2)
   :Y^\pm(z_1)\,Y^\pm(z_2)\,X^\pm(w):.\notag
\end{align}
} % end of \allowdisplaybreaks

%%%%%%%%%%%%%%%%%%%%%%%%%%%%%%%%%%%%%%%%%%%%%%%%%%%%%%%%%%%%
\section{Surjection from $\uqp(\spfh)$ to $U'$}

Let us restrict $\uqp(\spfh)$ to level $1$ representations
and denote it by $(\uqp(\spfh))_1$.
Recall from~\cite{MR95i:17011,MR88j:17020,MR99j:17021}
that under the identification between the two presentations
of $\uq(\spfh)$, we have
\begin{equation}\label{eq:n51}
t_0  \leftrightarrow \gamma (K_1^2 K_2)^{-1},
\end{equation}
so that
\begin{equation}
q^{2c} = t_0 t_1^2 t_2 \leftrightarrow \gamma,
\end{equation}
where $c = h_0 + h_1 + h_2$ is the canonical central element.
This shows that restricting to level $1$ is equivalent
to setting $\gamma = q^2$ in $\uq(\spfh)$.
This section is devoted to giving a surjection
from the quantum affine algebra $(\uqp(\spfh))_1$ to $U'$.
This would imply that the irreducible representation $\calV$
constructed for $U'$ is also an irreducible representation
for $\uqp(\spfh)$.

\begin{thm}\label{thm:51}
The following map defines a surjection from $(\uqp(\spfh))_1$
to $U'$.
\begin{alignat*}{2}
X_1^\pm(z) &\mapsto X^\pm(z), & \qquad
X_2^\pm(z) &\mapsto Y^\pm(z),\\
a_1(\pm k) &\mapsto a(\pm k), & \qquad
a_2(\pm k) &\mapsto - \frac{1}{[2]}\big( a(\pm k) + [2k]\, b(\pm k) \big),\\
K_1 &\mapsto K,&
K_2 &\mapsto (q^{2b(0)}K)^{-1}.
\end{alignat*}
\end{thm}

The rest of this section is devoted to proving this theorem.
Since the surjectivity is obvious, it suffices to show
that the image under this map of every defining relation for $\uqp(\spfh)$,
given in Subsection~\ref{ssec:22},
is also satisfied in $U'$.

The $i=j=1$ cases of these relations are immediate from the defining
relations of $\uq(\slth)$, so we shall not mention them below.
The first few equations are easy to check, so let us start
with~\eqref{eq:211}.

\vspace{2mm}\noindent\emph{Equation~\eqref{eq:211}} :
The $i=1$, $j=2$ case follows from~\eqref{eq:418}.
To check the $i=2$, $j=1$ case, we refer to~\eqref{eq:24} and
%it suffices to calculate
%\begin{equation*}
%- \frac{1}{[2]} \Big( \pm \frac{[2k]}{k} q^{\mp|k|}z^k \Big)
%= \pm \frac{[-k]_2}{k} q^{\mp|k|}z^k.
%\end{equation*}
notice that $-\frac{[2k]}{[2]} = [-k]_2$.
Last of all, to check the $i=j=2$ case, by use of~\eqref{eq:418}
and~\eqref{eq:419}, it suffice to calculate
\begin{align*}
-\frac{1}{[2]}\Big(
\mp\frac{[2k]}{k} q^{\mp |k|} &z^k \mp [|2k|] \frac{(q^{2k}-1+q^{-2k})}{|k|}
q^{\mp|k|}z^k \Big)\\
&= \pm \frac{1}{[2]}\frac{[2k]}{k} (q^{2k} + q^{-2k}) q^{\mp|k|}z^k\\
&= \pm \frac{[2k]_2}{k} q^{\mp|k|}z^k.
\end{align*}

\vspace{2mm}\noindent\emph{Equation~\eqref{eq:212}} :
This follows from the equations \eqref{eq:410}--\eqref{eq:417}.
For example, from~\eqref{eq:414},
one can show
\begin{equation*}
(z - q^4 w)\, Y^+(z) Y^+(w) + (w - q^4 z)\, Y^+(w) Y^+(z) = 0.
\end{equation*}

\vspace{2mm}\noindent\emph{Equation~\eqref{eq:213}} :
The $i=1$, $j=2$ case follows from~\eqref{eq:327}.
Likewise, the $i=2$, $j=1$ case follows from~\eqref{eq:336}.
To check the remaining $i=j=2$ case, we refer to~\eqref{eq:415}
and~\eqref{eq:416}.
We may calculate
\begin{align*}
\frac{1}{z^2}&\frac{1}{(1 - q^2 w/z)(1 - w/q^2 z)} -
\frac{1}{w^2}\frac{1}{(1 - q^2 z/w)(1 - z/q^2 w)}\\
&=
\delta(z /q^2 w)\frac{1}{w^2}\frac{w/q^2 z}{1-w/q^2 z} +
\delta(q^2 z/w)\frac{1}{w^2}\frac{q^2 w/z}{1-q^2 w/z}\\
&=
\frac{1}{q^2 - q^{-2}} \Big(
\delta(z /q^2 w) - \delta(q^2 z/w) \Big) \frac{1}{zw}
\end{align*}
and the rest follows if we carefully write down $:Y^+(z)Y^-(w):$\,.

\vspace{2mm}\noindent\emph{Equation~\eqref{eq:214}} :
This one is long.
Let us just sketch the $+$ part.
First, substitute the bosonization of $Y^+(w)$ in place of $X_2^+(w)$
and substitute $X^+(z_i)$ in place of $X_1^+(z_i)$ in~\eqref{eq:214}.
If we then send all $a(k)$ with $k < 0$ to the left of all the $X^+(z_i)$
and send all other parts of $Y^+(w)$ to the right
of all the $X^+(z_i)$, then we are left with showing
that
\begin{equation*}
\Sym_{z_1, z_2, z_3}
\begin{pmatrix}
A(z_1, z_2, z_3, w)
\hfill\\[0.2em]
\times\exp( - {\textstyle\sum} a_- w^+ )
\hfill\\[0.2em]
\times X^+(z_1)\, X^+(z_2)\, X^+(z_3)
\hfill\\[0.2em]
\times \exp( + {\textstyle\sum} a_+ w^- ) \Dmh (w)^{- h_1}
\hfill
\end{pmatrix}
= 0,
\end{equation*}
where
\begin{equation*}
A(z_1, z_2, z_3, w) =
\begin{pmatrix}
\hfill \displaystyle
\frac{z_1 z_2 z_3}{w^3}
 \frac{1}{(1-z_1/q^2 w)}\frac{1}{(1-z_2/q^2 w)}\frac{1}{(1-z_3/q^2 w)}\\[1em]
\hfill \displaystyle
+ [3] \frac{z_2 z_3}{w^2}
 \frac{1}{(1-w/q^2 z_1)}\frac{1}{(1-z_2/q^2 w)}\frac{1}{(1-z_3/q^2 w)}\\[1em]
\hfill \displaystyle
+ [3] \frac{z_3}{w}
 \frac{1}{(1-w/q^2 z_1)}\frac{1}{(1-w/q^2 z_2)}\frac{1}{(1-z_3/q^2 w)}\\[1em]
\hfill \displaystyle
+
 \frac{1}{(1-w/q^2 z_1)}\frac{1}{(1-w/q^2 z_2)}\frac{1}{(1-w/q^2 z_3)}
\end{pmatrix}.
\end{equation*}
Here, we have written the bosonized parts of $Y^+(w)$ in a simplified
manner.
Now, if we use the simple relation
\begin{equation*}
\frac{y/x}{1 - y/x} = \delta(x/y) - \frac{1}{1 - x/y},
\end{equation*}
we may show
\begin{equation*}
A(z_1, z_2, z_3, w) =
\frac{z_1 z_2 z_3}{w^3}\big(
B(z_1, z_2, z_3, w) + C(z_1, z_2, z_3, w)
\big),
\end{equation*}
with 
\begin{equation*}
B(z_1, z_2, z_3, w) =
\begin{pmatrix}
\hfill \displaystyle
 \frac{1}{(1-z_1/q^2 w)}\frac{1}{(1-z_2/q^2 w)}\frac{1}{(1-z_3/q^2 w)}\\[1em]
\hfill \displaystyle
- [3] q^2
 \frac{1}{(1-q^2 z_1/w)}\frac{1}{(1-z_2/q^2 w)}\frac{1}{(1-z_3/q^2 w)}\\[1em]
\hfill \displaystyle
+ [3] q^4
 \frac{1}{(1-q^2 z_1/w)}\frac{1}{(1-q^2 z_2/w)}\frac{1}{(1-z_3/q^2 w)}\\[1em]
\hfill \displaystyle
- q^6
 \frac{1}{(1-q^2 z_1/w)}\frac{1}{(1-q^2 z_2/w)}\frac{1}{(1-q^2 z_3/w)}
\end{pmatrix}
\end{equation*}
and
\begin{equation*}
C(z_1, z_2, z_3, w) =
\begin{pmatrix}
\hfill \displaystyle
[3] q^2
 \delta(q^2 z_1/w) \frac{1}{(1-z_2/q^2 w)}\frac{1}{(1-z_3/q^2 w)}\\[1em]
\hfill \displaystyle
+ [3] q^4
 \delta(q^2 z_1/w)
 \frac{w/q^2 z_2}{(1-w/q^2 z_2)}\frac{1}{(1-z_3/q^2 w)}\\[1em]
\hfill \displaystyle
- [3] q^4
 \frac{1}{(1-q^2 z_1/w)} \delta(q^2 z_2/w)\frac{1}{(1-z_3/q^2 w)}\\[1em]
\hfill \displaystyle
+ q^6
 \delta(q^2 z_1/w) \frac{w/q^2 z_2}{(1-w/q^2 z_2)}
  \frac{w/q^2 z_3}{(1-w/q^2 z_3)}\\[1em]
\hfill \displaystyle
- q^6
 \frac{1}{(1-q^2 z_1/w)}\delta(q^2 z_2/w)\frac{w/q^2 z_3}{(1-w/q^2 z_3)}\\[1em]
\hfill \displaystyle
+ q^6
 \frac{1}{(1-q^2 z_1/w)}\frac{1}{(1-q^2 z_2/w)} \delta(q^2 z_3/w)
\end{pmatrix}.
\end{equation*}
It suffices to show
\begin{equation}\label{eq:bzero}
\Sym_{z_1, z_2, z_3}\big(
B(z_1, z_2, z_3, w)\, X^+(z_1)\, X^+(z_2)\, X^+(z_3)
\big) = 0,
\end{equation}
and
\begin{equation}\label{eq:czero}
\Sym_{z_1, z_2, z_3}\big(
C(z_1, z_2, z_3, w)\, X^+(z_1)\, X^+(z_2)\, X^+(z_3)
\big) = 0.
\end{equation}
These may be done through extensive use of relation~\eqref{eq:25}.
For example, to do the first one, we first unify denominators.
There are many cancellation of terms and we have
\begin{equation*}
B(z_1, z_2, z_3, w) =
\begin{pmatrix}
(q^2 - q^{-2})\hfill\\[0.2em]
\displaystyle
\times\frac{1}{(1-z_1/q^2 w)}\frac{1}{(1-z_2/q^2 w)}\frac{1}{(1-z_3/q^2 w)}
\hfill\\[1em]
\displaystyle
\times\frac{1}{(1-q^2 z_1/w)}\frac{1}{(1-q^2 z_2/w)}\frac{1}{(1-q^2 z_3/w)}
\hfill\\[1em]
\displaystyle
\times
\begin{pmatrix}
\displaystyle
  \phantom{+}
  \frac{1}{w}\big( - {z_1} + (q^2 + q^4){z_2} - q^6{z_3} \big)
  \hfill\\[1em]
\displaystyle
+ \frac{1}{w^2}\big( - {z_1 z_2} - q^6 {z_2 z_3} + (q^2 + q^4) {z_1 z_3}\big)
  \hfill
\end{pmatrix}
\end{pmatrix}
.
\end{equation*}
We throw away all symmetric parts, and are left with showing
\begin{align}
&\Sym_{z_1, z_2, z_3}\big\{
(- {z_1} + (q^2 + q^4){z_2} - q^6{z_3})\,
X^+(z_1)\, X^+(z_2)\, X^+(z_3)
\big\} = 0,\label{eq:expan}\\
&\Sym_{z_1, z_2, z_3}\big\{
( - {z_1 z_2} - q^6 {z_2 z_3} + (q^2 + q^4) {z_1 z_3})\,
X^+(z_1)\, X^+(z_2)\, X^+(z_3)
\big\} = 0.\label{eq:expant}
\end{align}
The first symmetrization expands into
\begin{align*}
&( - z_1 + (q^2 + q^4) z_2 - q^6 z_3 )\, X^+(z_1)\, X^+(z_2)\, X^+(z_3)\\
+
&( - z_1 + (q^2 + q^4) z_3 - q^6 z_2 )\, X^+(z_1)\, X^+(z_3)\, X^+(z_2)\\
+
&( - z_3 + (q^2 + q^4) z_2 - q^6 z_1 )\, X^+(z_3)\, X^+(z_2)\, X^+(z_1)\\
+
&( - z_2 + (q^2 + q^4) z_1 - q^6 z_3 )\, X^+(z_2)\, X^+(z_1)\, X^+(z_3)\\
+
&( - z_2 + (q^2 + q^4) z_3 - q^6 z_1 )\, X^+(z_2)\, X^+(z_3)\, X^+(z_1)\\
+
&( - z_3 + (q^2 + q^4) z_1 - q^6 z_2 )\, X^+(z_3)\, X^+(z_1)\, X^+(z_2).
\end{align*}
In this long equation, we may single out
\begin{align*}
&( - z_1 + q^2 z_2 ) X^+(z_1)\, X^+(z_2)\, X^+(z_3) +
 ( - z_2 + q^2 z_1 )  X^+(z_2)\, X^+(z_1)\, X^+(z_3) = 0,\\
&( q^4 z_2 - q^6 z_3 ) X^+(z_1)\, X^+(z_2)\, X^+(z_3) +
 ( q^4 z_3 - q^6 z_2 ) X^+(z_1)\, X^+(z_3)\, X^+(z_2) = 0,\\
&( - z_1 + q^2 z_3 )  X^+(z_1)\, X^+(z_3)\, X^+(z_2) +
 ( - z_3 + q^2 z_1 ) X^+(z_3)\, X^+(z_1)\, X^+(z_2) = 0,
\end{align*}
and so on.
Note that these are all zero, thanks to~\eqref{eq:25}.
So, equation~\eqref{eq:expan} is indeed true.
The second symmetrization~\eqref{eq:expant}
may similarly be shown to be equal to zero,
leading to~\eqref{eq:bzero}.

Showing the validity of~\eqref{eq:czero} is quite similar, except that
we have the added complexity of using
\begin{equation*}
x\, \delta(x/y) = y\, \delta(x/y)
\end{equation*}
at appropriate places.

\vspace{2mm}\noindent\emph{Equation~\eqref{eq:215}} :
With equations~\eqref{eq:420}--\eqref{eq:422} in hand, it suffices to show
that the symmetrization of
\begin{equation*}
\begin{pmatrix}
\hfill\displaystyle
\frac{1}{(1 - q^{\mp 2}z_1/w)}\frac{1}{(1 - q^{\mp 2}z_2/w)}\\[1em]
\hfill\displaystyle
+ [2]_2 \frac{w}{z_1}
   \frac{1}{(1 - q^{\mp 2}w/z_1)}\frac{1}{(1 - q^{\mp 2}z_2/w)}\\[1em]
\hfill\displaystyle
+\frac{w^2}{z_1 z_2}
   \frac{1}{(1 - q^{\mp 2}w/z_1)}\frac{1}{(1 - q^{\mp 2}w/z_2)}
\end{pmatrix}
\times (z_1 - q^{\mp 4}z_2)(z_1 - z_2)
\end{equation*}
with respect to $z_1$ and $z_2$ is equal to zero.
It is done through routine calculation.

This concludes the proof of Theorem~\ref{thm:51}.

%%%%%%%%%%%%%%%%%%%%%%%%%%%%%%%%%%%%%%%%%%%%%%%%%%%%%%%%%%%%
\section{Irreducible highest weight modules}

We show in this section that the irreducible modules created
in Section~\ref{sec:4} are actually highest weight modules
of $\uq(\spfh)$.

Let us use the notation
\begin{equation*}\text{}
[a,b]_v = a b - v\, b a.
\end{equation*}
It is easy to check that
\begin{equation*}\text{}
[a,[b,c]_u]_v = [[a,b]_x, c]_{uv/x} + x [b,[a,c]_{v/x}]_{u/x},
\end{equation*}
for $x\neq 0$ (\cite{MR99j:17021}).

We first prepare a small lemma.
\begin{lem}\label{lem:61}
As operators acting on the highest weight vector $v_\la \in V(\la)$,
\begin{align*}
\text{}
[ x^-(0), [\Phitzz, &x^-(1)]_{q^{-2}}]_{1}\\
&=
- q^{-2} z^{-1} \Dh f_1^2
 \exp\bigg( \sum_{k=1}^{\infty}
            \frac{q^{k}}{[2k]} a(-k) z^k \bigg)
 z^{\half h_1}.
\end{align*}
\end{lem}
\begin{proof}
We may calculate
\begin{align*}
\text{}
[ x^-(0),\ &[\Phitzz, x^-(1)]_{q^{-2}}]_{1}\\
&=
[[ x^-(0), \Phitzz]_{q^{-2}}, x^-(1)]_{1}
+
q^{-2}[\Phitzz, [x^-(0), x^-(1)]_{q^2}]_{1}.
\end{align*}
Notice that the second term is zero because,
\begin{equation*}
\text{}
[x^-(0), x^-(1)]_{q^2}
= [ f_1, e_0 t_0^{-1}]_{q^2}
= [ f_1, e_0 ]_{1} t_0^{-1}
= 0.
\end{equation*}
And from~\eqref{eq:310}, we have
\begin{equation*}
\text{}
[ x^-(0), \Phitzz]_{q^{-2}}
=
q^{-2} z^{-1} [x^-(1), \Phitzz]_{q^2}
\end{equation*}
with $x^-(1) = e_0 t_0^{-1}$.
So, when acting on an extremal vector,
\begin{align*}
\text{}
[ x^-(0),\ &[\Phitzz, x^-(1)]_{q^{-2}}]_{1}\\
&= - q^{-2} z^{-1} x^-(1) x^-(1) \Phitzz,\\
&=
\begin{pmatrix}
\displaystyle
- q^{-2} z^{-1} x^-(1) x^-(1)\hfill\\[0.1em]
\displaystyle
   \times\exp\bigg( \sum_{k=1}^{\infty}
     \frac{q^{k}}{[2k]} a(-k)\, z^k \bigg)
   \exp\bigg( - \sum_{k=1}^{\infty}
     \frac{q^{k}}{[2k]} a(k)\, z^{-k} \bigg)\hfill\\[1.2em]
\displaystyle
   \times\Dh z^{\half h_1}\hfill
\end{pmatrix}
,\\
&= 
- q^{-2} z^{-1} x^-(1) x^-(1) \Dh
 \exp\bigg( \sum_{k=1}^{\infty}
            \frac{q^{k}}{[2k]} a(-k) z^k \bigg)
 z^{\half h_1},\\
&=
- q^{-2} z^{-1} \Dh x^-(0) x^-(0)
 \exp\bigg( \sum_{k=1}^{\infty} 
            \frac{q^{k}}{[2k]} a(-k) z^k \bigg)
 z^{\half h_1}.
\end{align*}
We have used Theorem~\ref{thm:34} and
equations~\eqref{eq:n309} and~\eqref{eq:n313t}
in the above line of equalities.
We've also used the fact $a(k) v_\la = 0$ for all $k > 0$.

Substitute $f_1 = x^-(0)$ to complete the proof.
\end{proof}

We now give a theorem that contains most of the result we have wanted.
\begin{thm}
For each $j=0,1,2$, the irreducible $\uqp(\spfh)$-module $\calV(j)$
is a highest weight module of highest weight weight $\La_j$.
The highest weight vectors are given as follows.
\begin{align*}
v_{2\La_0} \ot v(0)
  &\in \calV(0),\\
v_{\La_0 + \La_1} \ot v(- 1/2)
  &\in \calV(1),\\
v_{2\La_0} \ot v(-1)
  &\in \calV(2).
\end{align*}
\end{thm}
\begin{proof}
Using Theorem~\ref{thm:51} and
equation~\eqref{eq:n51}, it is straightforward
to check that the weight of each proposed highest weight vector
is $\La_j$.

It suffices to show that these are
indeed killed by each $e_i$ ($i=0,1,2$).
Recall from~\cite{MR95i:17011,MR88j:17020,MR99j:17021} that,
\begin{align*}
e_0 &= q^2\, [ x_1^-(0), [x_2^-(0), x_1^-(1)]_{q^{-2}}]_{1}
       K_1^{-2} K_2^{-1},\\
e_1 &= x_1^+(0),\\
e_2 &= x_2^+(0),
\end{align*}
for elements of $(\uq(\spfh))_1$.
Let us try the last one, $w = v_{2\La_0} \ot v(-1)$, as an example.
Other cases are simpler.

Action of $e_1$ on $w$ is trivially zero.
To see the action of $e_2$, we look for the coefficient of $z^{-1}$
in $Y^+(z)\cdot w$.
As $b(k) v(-1) = 0$ for all $k>0$,
we have
\begin{equation*}
\Omega_2(z) v(-1) =
 z \exp\bigg(-\sum_{k=1}^\infty b(-k) q^{-k} z^k\bigg) v(-2).
\end{equation*}
Hence, the smallest power of $z$ with nonzero coefficient appearing
from the $\Omega_2(z)$ part of $Y^+(z)\cdot w$ is $z$.
Likewise, we have
\begin{equation*}
\Psi^{(2)}_2(q^{-2}z) v_{2\La_0} =
  \Dmh
  \exp\bigg( - \sum_{k=1}^{\infty}
      \frac{q^{-k}}{[2k]} a(-k) z^k\bigg)
  v_{2\La_0}.
\end{equation*}
The smallest power possible is $z^{0}$.
Combined, they imply that the smallest power of $z$ appearing
in $Y^+(z) \cdot w$ is at least $z^1$.
The coefficient of $z^{-1}$ in $Y^+(z)\cdot w$
is zero and hence $e_2 w = 0$.

It remains to check the action of $e_0$.
Again, we look for the coefficient of $z^{-1}$ in
\begin{equation*}
q^2\, [x^-(0), [Y^-(z), x^-(1)]_{q^{-2}}]_{1}
       K_1^{-2} K_2^{-1} \cdot w.
\end{equation*}
As before, we may verify that
the smallest power of $z$ with nonzero coefficient appearing
from the $\Omega_0(z)$ part is $z^{-1}$.
Using Lemma~\ref{lem:61}, we see that what remains is
\begin{equation*}
- z^{-1} \Dh f_1^2 
\exp\bigg( \sum_{k=1}^{\infty}
            \frac{q^{k}}{[2k]} a(-k) z^k \bigg)
 z^{\half h_1}
K^{-1} v_{2\La_0}.
\end{equation*}
Since $f_1 v_{2\La_0} = 0$, the smallest nonzero power of $z$ is
(or, at least for the moment, seems to be)
$z^0$ with the coefficient
\begin{equation*}
- \Dh f_1^2 \frac{q}{[2]} a(-1)\, v_{2\La_0}.
\end{equation*}
We ignore the insignificant coefficients and follow
{\allowdisplaybreaks
\begin{align*}
\Dh f_1^2 a(-1)\, v_{2\La_0}
&= \Dh \big(x^-(0)\big)^2 a(-1)\, v_{2\La_0},\\
&= \big(x^-(1)\big)^2 a(-1) \Dh\, v_{2\La_0},\\
&= \big(x^-(1)\big)^2 a(-1)\, v_{2\La_1},\\
&= x^-(1) \big( a(-1) x^-(1) + q[2] x^-(0) \big)\, v_{2\La_1},\\
&= e_0 t_0^{-1} \big( a(-1) e_0 t_0^{-1} + q[2] f_1 \big)\, v_{2\La_1},\\
&= q[2]\, e_0 t_0^{-1} f_1\, v_{2\La_1},\\
&= q^3 [2]\, t_0^{-1} f_1 e_0\, v_{2\La_1},\\
&= 0.
\end{align*}}
Hence, the smallest nonzero power of $z$ was not $z^0$, but $z^1$.
Recalling that the smallest power obtained from $\Omega_0(z)$ part
was $z^{-1}$, we conclude that
the coefficient of $z^{-1}$, which we have been looking for, is zero.

We have $e_0 w = 0$
and the vector $w = v_{2\La_0} \ot v(-1)$ is an extremal vector.
This completes the proof.
\end{proof}

Now, note that $\uqp(\spfh)$ is a $d$-graded algebra.
And we know that the irreducible highest weight modules are
quotients of $\uqp(\spfh)$ by $d$-graded submodules.
So setting the degree of each highest weight vector to
zero determines the action of $q^d$ uniquely to the
rest of the irreducible module.

We have thus reached our goal.
\begin{thm}
$\calV(j)$ is the irreducible highest weight module
over $\uq(\spfh)$ of highest weight $\La_j$, for each $j=0,1,2$.
\end{thm}

We end this paper by writing down the explicit action of $q^d$ on
our new realization of highest weight modules.
First, add an element $q^d$ to the Heisenberg algebra $\Heis'$
with the relation
\begin{equation*}
q^d b(k) q^{-d} = q^k b(k)
\end{equation*}
The resulting algebra is denoted by $\Heis$.
We next make each $\calV_p$ into a $\Heis$-module.
It suffices to define the action of $q^d$ on each $v(p)$.
This is done differently, depending on the irreducible
highest weight module we want to create.
\begin{alignat*}{2}
&\calV(0):\qquad &
q^d v(p) &=
\begin{cases}
q^{-\half p^2} v(p) & \text{for $p\in 2\Z$,}\\
q^{-\half(p^2 + 1)} v(p) & \text{for $p\in 2\Z + 1$.}
\end{cases}\\
&\calV(1):\qquad &
q^d v(p) &=
q^{-\half(p-\half)(p+\half)} v(p) \quad\text{for $p\in \Z + \half$.}\\
&\calV(2):\qquad &
q^d v(p) &=
\begin{cases}
q^{-\half p^2} v(p) & \text{for $p\in 2\Z$,}\\
q^{-\half(p^2 - 1)} v(p) & \text{for $p\in 2\Z + 1$.}
\end{cases}
\end{alignat*}
Finally, the action of $q^d$ on $\calV(j)$ ($j=0,1,2$) is
then equal to
\begin{equation*}
q^d \ot q^d \in \End{\calV(j)}.
\end{equation*}
Here, $q^d$ on the left of the tensor sign is the usual
$\uq(\slth)$-action, and $q^d$ on the right of the tensor sign
signifies the action just defined.

%%%%%%%%%%%%%%%%%%%%%%%%%%%%%%%%%%%%%%%%%%%%%%%%%%%%%%%%%%%%
\providecommand{\bysame}{\leavevmode\hbox to3em{\hrulefill}\thinspace}

\end{document}